\documentclass[a4paper,11pt,makeidx]{amsart}
\usepackage{amscd}
\usepackage{xypic}  %commu
\usepackage{amssymb}
\usepackage{amsthm}
\usepackage{epsf}
\makeindex

\newtheorem{thm}{Theorem}[section]
\newtheorem{cor}[thm]{Corollary}
\newtheorem{exa}[thm]{Example}
\newtheorem{lemma}[thm]{Lemma}

\newtheorem{prop}[thm]{Proposition}
\newtheorem{lemma-exa}[thm]{Lemma-Example}

\def\Hilb{\operatorname{Hilb}}
\def\NS{\operatorname{NS}}
\def\kod{\operatorname{kod}}
\def\Sing{\operatorname{Sing}}

\def\im{\operatorname{im}}

\def\length{\operatorname{length}}
\def\c1{\operatorname{c_1}}
\def\c2{\operatorname{c_2}}

\def\Rat{\operatorname{Rat}}

\def\Sym{\operatorname{Sym}}

\def\Alb{\operatorname{Alb}}

\def\ZZ{{\mathbb Z}}
\def\NN{{\mathbb N}}

\def\PP{{\mathbb P}}
\def\FF{{\mathbb F}}

\def\C{{\mathcal C}}

\def\N{{\mathcal N}}
\def\O{{\mathcal O}}

\def\T{{\mathcal T}}

\def\*{\otimes}
\def\x{\times}                  % product (fiber)
\def\iso{\simeq}
\def\eqv{\equiv}
\def\sub{\subseteq}

\def\sup{\supseteq}

\def\+{\oplus}                   % direct sum
\def\*{\otimes}                  % tensor product
\def\hpil{\longrightarrow}       % ----->
\def\khpil{\rightarrow}

\def\Pic{\operatorname{Pic}}

\def\Supp{\operatorname{Supp}}

\begin{document}

\title[Families of curves with hyperelliptic normalizations]{Remarks on families of singular curves with hyperelliptic normalizations}
\author{Andreas Leopold Knutsen}

\address{\hskip -.43cm Andreas Leopold Knutsen, Dipartimento di Matematica, Universit\`a di Roma Tre, Largo San Leonardo
Murialdo 1, 00146, Roma, Italy. e-mail {\tt knutsen@mat.uniroma3.it}}

\thanks{Research supported by a Marie Curie Intra-European Fellowship within the 6th European
Community Framework Programme}

\thanks{{\it 2000 Mathematics Subject Classification}: Primary 14H10, Secondary 14E30, 14H51, 14J10}

\begin{abstract}
 We give restrictions on the existence of families of 
curves on smooth projective surfaces $S$ of nonnegative Kodaira dimension all having constant geometric genus
$g \geq 2$ and hyperelliptic normalizations. In particular, we
prove a Reider-like result whose proof is ``vector bundle-free'' and relies on deformation theory and bending-and-breaking of rational curves in $\Sym^2(S)$. We also give examples of families of such curves.
\end{abstract}

\maketitle

\section{Introduction} \label{sec:intro}

The object of study of this  paper is families of irreducible curves
with hyperelliptic normalizations (of genus $\geq 2$) on a smooth surface $S$. Such families give rise to, 
because of the unique $g^1_2$s on their normalizations, families of the same dimensions of irreducible rational curves in the Hilbert scheme $\Hilb^2(S)$. Because of the importance of rational curves and the subvarieties they cover due to Mori theory, it is natural to try to check the existence of, or bound the dimensions of, families of such
curves, or alternatively, their counterparts on $S$.

Let  $S$ be a smooth surface and $V$ be a reduced and irreducible scheme parametrizing a flat family of curves on $S$ all having constant geometric genus
$g \geq 2$ and hyperelliptic normalizations. 

It is easy to see (cf. Lemma \ref{lemma:0} below) that
if $|K_S|$ is birational, then $\dim V=0$. This shows that the problem of bounding the dimension of a family of curves with hyperelliptic normalizations is solved for a large class of surfaces. At the same time,
it is relatively easy to find obvious examples of surfaces with large families of
curves with hyperelliptic normalizations: 
In fact, if $S$ is any smooth surface admitting a rational $2:1$ map $f:S \rightarrow R$ onto a rational surface, then we can just pull back families of rational curves on $R$. 
There are several such examples of such double covers, even for $p_g(S) >0$, 
see for instance the works of Horikawa \cite{Hor1}-\cite{Hor4} for surfaces of general type, Saint-Donat  \cite{SD} for $K3$ surfaces, and \cite{se,sv,BFL} 
for surfaces with smooth hyperelliptic hyperplane sections.

We note that {\it smooth} hyperelliptic curves on surfaces have been extensively studied by means of adjunction theory (see \cite{se,sv,BFL} to mention a few). Of course Reider's famous result \cite{re}
can be used to prove that if $C \subset  S$ is a smooth hyperelliptic curve and $C^2 \geq 9$, then 
there is a pencil $|E|$ such that either $E^2=0$ and $E.C=2$, or $E^2=1$, $C \eqv 3E$ and $|E|$ has one base point $x$ lying on $C$ (making the obvious modifications in the proof of 
\cite[Cor. 2]{re}). Unfortunately no such results seem to be available, at least as far as we know, in the case of singular curves.

In this paper we prove some results bounding dimensions of families of irreducible curves
with hyperelliptic normalizations on smooth surfaces $S$ in Section \ref{sec:exa}. In particular, we show that
the dimension is bounded by one if $S$ is fibered over a smooth nonhyperelliptic curve of genus $\geq 3$
(Lemma \ref{lemma:sucurva}) and by two if $S$ has maximal albanese dimension (Proposition \ref{prop:dimboundalb}). We also give several examples of families of irreducible curves
with hyperelliptic normalizations. Then, in Section \ref{sec:reider}, we prove a Reider-like result, cf. Theorem \ref{thm:reider},
stating that
any family of dimension $\geq 3$ (resp. $\geq 5$) of curves with hyperelliptic normalizations on a smooth surface $S$ with $p_g(S)>0$ (resp. $\kod(S) \geq 0$ and $p_g(S)=0$) forces the existence of some special divisors enjoying some particular intersection properties. Moreover, these divisors ``cut out'' the $g^1_2$s on the normalizations of the curves in the family. 

We hope the results will find more applications and also hope that the reader will find the method of proof of interest: in fact, unlike the results on smooth curves, which use adjunction theory and/or vector bundle methods, our method uses deformations of curves and bending and breaking of rational curves in $\Sym^2(S)$, a method we also used in \cite{fkp2}. Thus, the special divisor occurring in  Theorem \ref{thm:reider} is obtained as a component of a degenerated member at the border of the family. In this sense our method is perhaps more geometric and intuitive than Reider's method. 

In Section \ref{sec:exa} we state the general setting, show that the dimension of families of curves with hyperelliptic normalizations can be bounded in various cases and give some examples of such families.

Then, in Section \ref{sec:reider} we prove the Reider-like result, Theorem \ref{thm:reider}, 
passing by Proposition \ref{prop:algdec} and Lemma \ref{lemma:algdec}. Finally, we make some remarks in 
Section \ref{sec:appexa},
including writing out the results in the case of smooth curves in \S\;\ref{ss:d0} recovering Reider's result, and in the case of only one singular point of multiplicity one in
\S\;\ref{ss:d1}.

\vspace{0,3cm}

\noindent {\it Acknowledgements.} I thank F.~Flamini, T.~Johnsen and L.~Stoppino for valuable suggestions and M.~Beltrametti for reading through a preliminary version of this paper. I also thank C.~Ciliberto for showing me a proof of what is now Lemma \ref{lemma:dimbound}.

\section{Dimension bounds and examples} \label{sec:exa}

Consider the following assumptions:
\begin{eqnarray} 
\label{eq:ass1} & \mbox{$S$ is a smooth projective surface with $\kod(S) \geq 0$ and $V$ a reduced 
and irreducible} & \\
\nonumber & \mbox{scheme parametrizing a flat family of irreducible curves on $S$ of constant arithmetic } & \\
\nonumber & \mbox{and geometric genera $p_a$ and $p_g \geq 2$, respectively, and hyperelliptic normalizations.} & \\
\nonumber & \mbox{We denote by $C$ the algebraic equivalence class of the curves.} &
\end{eqnarray}

Note that $V$ as in \eqref{eq:ass1} is, by default, nonempty.

We have the following elementary result, already mentioned in the introduction:

\begin{lemma} \label{lemma:0}
Under the assumptions  \eqref{eq:ass1}, if $|K_S|$ is birational, then $\dim V=0$.
\end{lemma}

\begin{proof}
We may assume that $\dim V=1$. Then after compactifying and resolving the singularities of 
the universal family over $V$, we obtain a smooth surface $T$, fibered over a smooth curve, with general fiber 
$F$ a smooth hyperelliptic curve of genus $\geq 2$, and a surjective morphism $f:T \to S$. 
By adjunction $|K_T|$ is not birational on the general fiber $F$. Since $K_T = f^*K_S + R$, where $R$ is the (effective) ramification divisor of $f$, and
$f$ is generically $1:1$ on the fibers,
we see that $|K_S|$ cannot be birational.
\end{proof}

Note that, as mentioned in the introduction, any irreducible curve $C$ on a surface $S$
with hyperelliptic normalization (of geometric genus $\geq 2$) gives rise to a unique irreducible rational curve $R_C \subset  \Hilb^2(S)$. Precisely, this can be seen in the following way: Let $\nu: \widetilde{C} \khpil C$ be the normalization. Then  the unique $g^1_2$ on $\widetilde{C}$ 
induces a $\PP^1 \subset  \Sym^2(\widetilde{C})$ and this is mapped to an irreducible rational curve $\Gamma_C$ in $\Sym^2(S)$ by the natural composed morphism 
\[
\xymatrix{
\Sym^2(\widetilde{C}_X) \ar[r]^{\tilde{\nu}^{(2)}} & \Sym^2(C_X)  \ar@{^{(}->}[r] & \Sym^2(S).
}
\]
The irreducible rational curve $R_C \subset  \Hilb^2(S)$ is the strict transform by the Hilbert-Chow morphism
$\mu: \Hilb^2(S) \khpil \Sym^2(S)$ of this curve. Note that the Hilbert-Chow morphism  resolves $\Sing (\Sym^2(S)) \iso S$ and gives an obvious one-to-one correspondence between irreducible curves in $\Hilb^2(S)$ not contained in the exceptional locus (which is a $\PP^1$-bundle over $S$) and irreducible curves in $\Sym^2(S)$ not contained in $\Sing (\Sym^2(S))$.

The correspondence in the opposite direction, that is, from irreducible rational curves in $\Hilb^2(S)$ to curves in $S$ is more delicate and we refer to \cite[\S\;2]{fkp2} for details. Suffice it to say that
irreducible rational curves in $\Hilb^2(S)$ not contained in the exceptional locus 
give rise to curves on $S$ with rational, elliptic or hyperelliptic normalizations, by taking the 
(one-dimensional component of the) union of the supports of the points of the curve in $\Hilb^2(S)$ when we consider 
these as length-two schemes on $S$.

To extend the correspondence into families, we proceed as follows (cf. also \cite{fkp2}).

Given \eqref{eq:ass1}, let $\varphi:\C \khpil V$ be the universal family. 
Normalizing $\C$ we obtain, possibly restricting to an open dense subscheme
of $V$, a flat family $\tilde{\varphi} : \widetilde{\C}\to V$ of smooth hyperelliptic 
curves of genus $p_g \geq 2$, cf. \cite[Thm. 1.3.2]{tes}. 
Let $\omega_{\widetilde{\C}/V}$ be the relative dualizing sheaf. 
The morphism 
$\gamma: \widetilde{\C}\to \PP(\tilde{\varphi}_* (\omega_{\widetilde{\C}/V}))$ over $V$
is finite and of relative degree two onto its image (cf. also \cite[Thm. 5.5 (iv)]{LK}), which we denote by $\mathcal{P}_V$.
We now have a universal family $\psi: \mathcal{P}_V \khpil V$ 
of rational curves and since the points in these curves correspond to couples of points of $S$, possibly coinciding, we have a natural morphism $\Phi_V: \mathcal{P}_V \khpil \Sym^2(S)$. We define
\begin{equation} \label{eq:defrv}
 R_V:= \overline{ \im \Phi_V}  \; \; \; \; \; \mbox{(the Zariski closure)}.
\end{equation}
We have
\begin{equation} \label{eq:univ}
\xymatrix{ & \widetilde{\C} \ar[dl]_{\pi} \ar[dr]_{\tilde{\varphi}} \ar[r]^{\gamma}  & 
\mathcal{P}_V \ar[d]^{\psi} \ar[r]^{\Phi_V} & R_V & \hspace{-1cm} \subset  \Sym^2(S)\\
S  & & V, & & 
}
\end{equation}
where $\pi$ is the natural morphism. Note that 
the natural map $\mathcal{P}_V \hpil \Rat(\Sym^2(S))$ defined by $\Phi_V$ 
is generically finite (where $\Rat(\Sym^2(S))$ is the union of the components of 
$\Hilb(\Sym^2(S))$ whose general points correspond to reduced connected curves with 
rational components) and $\Phi_V$ 
maps no curve in the family
to $\Sing(\Sym^2(S)) \iso S$ by construction. Also note that $\dim R_V \leq 3$,
as $\Sym^2(S)$ is not uniruled, since $\kod(S) \geq 0$ (see e.g. \cite[Prop. 2.1]{HT}).

Given $V$ as in \eqref{eq:ass1}, we will call $\mu_*^{-1}R_V \subset \Hilb^2(S)$, the strict transform of $R_V$ by the Hilbert-Chow morphism, the {\it locus covered by the associated rational curves in $Hilb^2(S)$}.

We will make use of the following consequence of Mumford's well-known theorem on $0$-cycles on surfaces
\cite[Cor. p.~203]{Mum}, as generalized in \cite[Cor.~3.2]{fkp2}:

\begin{prop} \label{prop:mumford}
Assume that $S$ is a smooth surface with $p_g(S) >0$ and $R \subset \Sym^2(S)$ is a subvariety that is covered by a family
of rational curves of dimension $\geq 3$.

Then $R$ is a surface with rational desingularization.
\end{prop}

\begin{proof}
This is \cite[Prop. 3.6]{fkp2}.
\end{proof}

We will now give some results bounding the dimension of $V$ as in \eqref{eq:ass1} in various situations.

\begin{lemma} \label{lemma:sucurva}
 Assume \eqref{eq:ass1} and that $f:S \khpil B$ is a fibration over a nonhyperelliptic smooth curve $B$ of genus $\geq 3$.

Then $\dim V \leq 1$ with equality holding if and only if the general fiber of $f$ is a (smooth) hyperelliptic curve and $V$ parametrizes a subset of the fibers.
\end{lemma}

\begin{proof}
Let $\{C_v\}_{v \in V}$ be the family on $S$ given by $V$ and  $\{\Gamma_v=\Gamma_{C_v}\}_{v \in V}$ be the associated family of rational curves in $\Sym^2(S)$, given by $\psi$ and 
$\Phi_V$ as in \eqref{eq:univ}. 

Consider $f^{(2)}:\Sym^2(S) \khpil \Sym^2(B)$. Since $B$ is neither elliptic nor hyperelliptic, 
$\Sym^2(B)$ does not contain rational curves. Therefore, $f^{(2)}$ must contract every $\Gamma_v$ to a point, 
say $b_v+b'_v \in \Sym^2(B)$, with $b_v, b'_v \in B$. Letting $F_v:=f^{-1}b_v$ and
$F'_v:=f^{-1}b'_v$ denote the two fibers of $f$, we get that  $\Gamma_v$ is contained in the surface
$F_v+F'_v \subset \Sym^2(S)$. Hence $\Supp(C_v)= F_v \cup F'_v \subset S$. Since each $C_v$ is irreducible, the result follows.
\end{proof}

The following example shows that  the condition that $B$ is neither elliptic nor hyperelliptic is in fact necessary:

\begin{exa} \label{exa:baseiper}
{\rm Start with a Hirzebruch surface $\FF_e$, with $e \geq 0$, with $\Pic \FF_e \iso \ZZ[\Sigma] \+ \ZZ[F]$, where $\Sigma^2=-e$, $F^2=0$ and $\Sigma.F=1$. Choose integers $\alpha$ and $\beta$ such that
\begin{equation} \label{eq:coeff}
 \alpha \geq 2 \; \mbox{and} \; \beta \geq\{2, \alpha e+1\}.
\end{equation}
and take a general pencil in $|\alpha\Sigma + \beta F|$. Note that the conditions \eqref{eq:coeff} guarantee that the general element in the general such pencil is in fact irreducible. Now
take the blow up $\pi: \widetilde{\FF}_e \khpil \FF_e$
along the
$(\alpha\Sigma + \beta F)^2=\alpha(2\beta- e \alpha)$ base points of the pencil
and denote the exceptional curves by $E_i$. Then
\[ K_{\widetilde{\FF}_e} \sim -2\pi^*\Sigma-(e+2)\pi^*F + \Sigma E_i. \]
Set
\[ \widetilde{D}:= \pi^*(\alpha\Sigma + \beta F) -\Sigma E_i. \]
Then $\widetilde{D}^2=0$ and $|\widetilde{D}|$ is a pencil defining a fibration 
$g:\widetilde{\FF}_e \khpil \PP_1$. For any integer $l \geq 2$, choose a general 
$\Delta_l \in |\O_{\PP^1}(2l)|$ consisting of distinct points and let $R_l \in |2l\widetilde{D}|$ be the corresponding (smooth) divisor. Then
$\Delta_l$ and $R_l$ define two double covers $\nu$ and $\mu$, respectively, that are compatible, in the sense
that we have a commutative diagram:
\begin{equation} \label{eq:commuta}
\xymatrix{ 
T \ar[d]_f \ar[r]^{\mu} & \widetilde{\FF}_e \ar[d]^g \\
           B \ar[r]^{\nu} & \PP^1,
}
\end{equation}
where $T$ is a smooth surface, $B$ is a smooth curve and $f$ is induced by $g$ in the natural way.

By Riemann-Hurwitz, the genus of $B$ is $g(B) =l-1 \geq 1$ and $B$ is either elliptic or hyperelliptic.
As
\[ K_T \sim \mu^*(K_{\widetilde{\FF}_e}+l\widetilde{D}) = \mu^*((l-1)\widetilde{D}+\pi^*((\alpha-2)\Sigma+(\beta-e-2)F), \]
the conditions \eqref{eq:coeff} imply $\kod(S) \geq 0$, in fact even $p_g(S) >0$. 

The surface $\widetilde{\FF}_e$, being rational, contains families of irreducible (smooth) rational curves of arbitrarily high dimensions. Pulling them back on $T$ yields families of (smooth) hyperelliptic curves of arbitrarily high dimensions.}
\end{exa}

\vspace{0,3cm}

Note that the surfaces $T$ in the example have irregularity $q(T)=h^1(\O_T)=g(B)=l-1 \geq 1$. Thus the surfaces have arbitrarily high irregularity. Nevertheless, using the albanese map (cf. e.g. \cite[pp.~80-88]{be0}), we can prove the following bounds.

\begin{prop} \label{prop:dimboundalb}
Assume \eqref{eq:ass1} with $q(S):=h^1(\O_S) \geq 2$ and let $\alpha: S \khpil \Alb S$ be the albanese map.

If $\im \alpha$ is a curve $B$, then $\dim V \leq 1$ unless $B$ is hyperelliptic.

If $\im \alpha$ is a surface (i.e., $S$ is of {\it maximal albanese dimension}), then $\dim V \leq 2$.
\end{prop}

\begin{proof}
If $\im \alpha$ is a curve $B$, then $B$ is necessarily smooth of genus $q(S) \geq 2$ (cf. e.g. \cite[Prop. V.15]{be0}). Then we apply Lemma \ref{lemma:sucurva}.

If $\im \alpha$ is a surface $T$, we must have $p_g(S) >0$ by \cite[Lemme V.18]{be0}. 
Assume now, to get a contradiction, that $\dim V \geq 3$. Then $R_V$, defined in \eqref{eq:defrv}, is a surface with rational desingularization, by Proposition \ref{prop:mumford}.

As above, let $\{C_v\}_{v \in V}$ be the family on $S$ given by $V$ and  $\{\Gamma_v=\Gamma_{C_v}\}_{v \in V}$ be the associated family of rational curves in $\Sym^2(S)$, given by $\psi$ and 
$\Phi_V$ as in \eqref{eq:univ}. 

Consider $\alpha^{(2)}:\Sym^2(S) \khpil \Sym^2(T) \subset \Sym^2(\Alb S)$. As $\alpha$ does not contract $C_v$, for general 
$v \in V$, the surface $R_V$ is mapped by $\alpha^{(2)}$ to a surface $R'_V \subset \Sym^2(T)$. Let
$\Sigma: \Sym^2(\Alb S) \khpil \Alb S$ be the summation morphism. As $\Alb S$, being abelian, cannot contain rational curves, cf. e.g. \cite[Prop. 4.9.5]{bl}, each rational curve $\Gamma \subset R'_V \subset \Sym^2(T)$ must be contracted to a point by $\Sigma_{|R'_V}$, say $\Sigma(\Gamma)=p_{\Gamma} \in \Alb S$. Now all fibers
$\Sigma^{-1}p$, for $p \in \Alb S$, are isomorphic to the Kummer variety of $\Alb S$, cf. e.g. \cite[\S4.8]{bl} for the definition. As rational curves on Kummer varieties cannot move, by \cite[Thm. 1]{pi}, we must have that any family of rational curves on $R'_V$ has dimension $\leq \dim \Sigma(R'_V) \leq 2$.
But this contradicts the fact that $R_V$ has rational desingularization.

Therefore $\dim V \leq 2$, as desired.
\end{proof}

The following result shows that equality $\dim V=2$ is in fact attained on abelian surfaces, which have maximal albanese dimension.

\begin{lemma-exa} \label{lemma:abe-biell}
 Assume \eqref{eq:ass1} with $S$ abelian. 

Then $\dim V=2$ and the locus covered by the associated rational curves in $\Hilb^2(S)$ is a threefold birational to a $\PP^1$-bundle over $S$.

Furthermore, such families exist if and only if $S$ is simple (i.e., not the product of two elliptic curves).

In particular, a bielliptic surface does not contain families as in \eqref{eq:ass1}.
\end{lemma-exa}

\begin{proof}
Assume that $S$ is abelian. Consider the natural composed morphism
\[ 
\xymatrix{ 
\alpha: \Hilb^2(S) \ar[r]^{\mu} & \Sym^2(S) \ar[r]^{\hspace{0,4cm}\Sigma} & S,  
}
\]
where $\mu$ is the Hilbert-Chow morphism and $\Sigma$ is the summation morphism.
As above, let $\{\Gamma_v=\Gamma_{C_v}\}_{v \in V}$ be the associated family of rational curves in 
$\Sym^2(S)$ given by $V$.
As $S$ cannot contain rational curves
 (see e.g. \cite[Prop. 4.9.5]{bl}),
$\Sigma$ must contract each $\Gamma_v$ to a point $p_v$. Therefore, the strict 
transform
$R_v=R_{C_v} :=\mu_*^{-1}(\Gamma_v) \subset  \Hilb^2(S)$ is contained in the surface
$\alpha^{-1}(p_v) \subset  \Hilb^2(S)$. Now all such fibers of $\alpha$ over points in $S$ 
are isomorphic to the (desingularized) Kummer surface of $S$ 
(cf. \cite[\S\;7]{Be}, \cite{Be1} or
\cite[2.3]{H1} and e.g. \cite[\S\;10.2]{bl} and \cite[V.16]{BPV} for the definition).  
Since a Kummer surface is $K3$, rational curves do not move inside it (this also follows from Lemma \ref{lemma:dimbound} below). Therefore, the family is given by 
$\{R_v\}_{p \in S}$,
proving the first assertion.

If $S$ is simple, then it contains irreducible curves of geometric genus two, see e.g. \cite[Cor. 2.2]{ls}.

If $S = E_1 \x E_2$, with each $E_i$ a smooth elliptic curve, then $\Sym^2(S)=\Sym^2(E_1) \x \Sym^2(E_2)$
and each $\Sym^2(E_i)$ is an elliptic ruled surface. Any rational curve in $\Sym^2(S)$ not lying in
$\Sing(\Sym^2(S))$ is therefore mapped by projection to a rational curve in either 
$\Sym^2(E_1)$ or $\Sym^2(E_2)$, which has to be a fiber of the ruling. Therefore, the rational curve in
$\Sym^2(S)$ corresponds to a $g^1_2$ on one of the elliptic fibers of $S$, proving that there is no 
irreducible curve on $S$ with hyperelliptic normalization of geometric genus $\geq 2$.

If $S$ is bielliptic, there is a finite morphism $f:T \khpil S$ where $T$ is a product of two elliptic curves, cf. e.g. \cite[D{\'e}f. VI.19]{be0} or \cite[p.~199]{BPV}, whence $T$ is abelian. Clearly, $f$ is unramified, as 
$K_T \sim \O_T$, whence so is $f^{(2)}: \Sym^2(T) \khpil \Sym^2(S)$. Therefore, the family of rational curves 
 $\psi: \mathcal{P}_V \khpil V$ as in \eqref{eq:univ} is pulled back, via $f^{(2)}$, to two copies of the family in $\Sym^2(T)$. By what we proved above, the corresponding families of curves on $T$ consist of elliptic 
curves, whence the same holds on $S$.
\end{proof}

We conclude this section by giving some examples of families as in \eqref{eq:ass1} of high dimensions.

\begin{exa} \label{exa:p2}
 {\rm Let $W(n) \subset  |\O_{\PP^2}(n)|$ denote the} Severi variety {\rm of nodal, irreducible rational curves in $|\O_{\PP^2}(n)|$. Then $W(n)$ is irreducible of dimension $3n-1$, by a well-known result of Severi and Harris, cf. \cite[Thm. 1.1]{CS} and \cite{har}. For any integer $b \geq 3$, take a general smooth $B \in |\O_{\PP^2}(2b)|$,
so that, for any $n$, the general curve in $W(n)$ intersects $B$ transversally. 

Let $f:S \khpil \PP^2$ be the double cover defined by $B$, so that $S$ is a smooth surface and $f$ is branched along $B$. Setting $H:=f^*\O_{\PP^2}(1)$, we have $K_S \sim (b-3)H$. Let $V(n) \subset  |nH|$ be the subscheme parametrizing the inverse images of the curves in $W(n)$ that intersect $B$ transversally and 
let $p_a(n)=p_a(nH)$ and $p_g(n)$ denote the arithmetic and geometric genera of the curves in $V(n)$. Then
$V(n)$ satisfies the conditions in \eqref{eq:ass1} and 
\[ \dim V(n)=3n-1, \; p_a(n)=n(n+b-3)+1 \; \mbox{and} \; p_g(n)=bn-1, \]
the second equality by adjunction and the third by Riemann-Hurwitz.
Note that the elements of $V(n)$ have $\delta(n):=p_a(n)-p_g(n)=n(n-3)+2$ nodes (two over each of the nodes 
of the corresponding curves in $W(n)$.)

Also note that curves in $W(n)$ that are tangent to $B$ yield families in $|nH|$ with lower
geometric genera and lower dimensions.

Of course $\Hilb^2(S)$ contains a copy of the $\PP^2$, which is precisely the locus in $\Hilb^2(S)$ covered by the rational curves associated to the curves in $V(n)$. 
}
\end{exa}

\begin{exa} \label{exa:enr}
 {\rm Let $S$ be an Enriques surface. Then $S$ contains several elliptic pencils and we can always pick (at least) two such, $|2E_1|$ and $|2E_2|$, with $E_1.E_2=1$, cf. \cite[Thms. 3 and 3.2]{cos}. (By adjunction $E_1^2=E_2^2=0$ and it is well-known that $2E_i$ and $2E_i'$, where $E_i'$ denotes the unique element of $|E_i+K_S|$, are the two multiple fibers of the elliptic pencils.) 
Consider $H:=3E_1+E_2$; then $H^2=6$ and the base scheme of $|H|$ 
consists of two distinct points $x$ and $y$, where
\[x=E_1 \cap E_2' \; \mbox{and} \; y=E_1' \cap E_2 \]
(see \cite[Prop. 3.1.6 and Thm. 4.4.1]{cd}). Let $f: \widetilde{S} \khpil S$ be the blow up along
$x$ and $y$ and $E_x$ and $E_{y}$ the two exceptional divisors. Set $L:=f^*H-E_x-E_{y}$.
Then $L^2=4$ and $|L|$ is base point free and, by \cite[Thm. 4.5.2]{cd}, defines a morphism of degree two $\varphi_L: \widetilde{S} \khpil Q$
onto a smooth quadric $Q \subset  \PP^3$, which can be seen as the embedding of $\FF_0 \iso \PP^1 \x \PP^1$ by the complete linear system $|\ell_1+ \ell_2|$, where $\ell_1$ and $\ell_2$ are the two rulings. In particular, by construction, 
$L \sim \varphi_L^*(\O_Q(1)) \sim  \varphi_L^*(\ell_1+\ell_2)$. Furthermore, by 
\cite[Rem. 4.5.1 and Thm. 4.5.2]{cd} the pencil $|f^*(2E_1)|$ on $\widetilde{S}$ is mapped by 
$\varphi_L$ to $|\ell_1|$, so that $\varphi_L^*\ell_1 \sim f^*(2E_1)$. We therefore have
$\varphi_L^*\ell_2 \sim L - \varphi_L^*\ell_1 \sim f^*(E_1+E_2)-E_x-E_y$.

It follows that, for any $n \geq 1$, the general smooth rational curve
in $|\ell_1+n\ell_2| \iso \PP^{2n+1}$ yields by pullback by $\varphi_L$ a 
smooth hyperelliptic curve in $|f^*((n+2)E_1+nE_2)-n(E_x+E_y)|$ on $\widetilde{S}$ of genus
$3n+1$ by adjunction 
(or by Riemann-Hurwitz and the description of the ramification in \cite[Thm. 4.5.2]{cd}). 

Pushing down to $S$ we thus obtain subschemes, for each $n \in \NN$,
\[ V(n) \subset  |(n+2)E_1+nE_2|, \; \mbox{such that} \; \dim V(n)=2n+1, \]
parametrizing irreducible curves with hyperelliptic normalizations of geometric genera $p_g(n)$
 and arithmetic genera $p_a(n)$, where
\[ p_g(n)= 3n+1 \; \mbox{and} \; p_a(n)=\frac{1}{2}((n+2)E_1+nE_2)^2+1=n(n+2)+1. \]
Note that for each $n \geq 2$ all the curves in the family have precisely two singular points, located at $x$ and $y$, both of multiplicity $n$.

Of course $\Hilb^2(S)$ contains a rational surface birational to $\PP^1 \x \PP^1$, which is precisely the locus in $\Hilb^2(S)$ covered by the rational curves associated to the curves in $V(n)$.

One can repeat the construction with $E_1$ and $E_2$ interchanged, or with other elliptic pencils on the surface. Moreover, choosing smooth rational curves that are tangent to the branch divisor of $\varphi_L$, we can obtain families with lower geometric genera, that is, 
with more singularities.

Moreover, one can also repeat the process for $H:=2E_1+E_2$, which defines a rational $2:1$ map onto $\PP^2$, following the lines of the previous example. Note that the smooth curves in $|2E_1+E_2|$ form a $2$-dimensional family of} smooth {\rm hyperelliptic curves on $S$, by \cite[Cor. 4.5.1]{cd}.
}
\end{exa}

\begin{exa} \label{exa:k3}
 {\rm Let $S$ be a $K3$ surface. Then $S^{[2]}$ is a} hyperk\"ahler
fourfold, {\rm also called an} irreducible symplectic fourfold, {\rm and rational curves and uniruled subvarieties are central in the study of the (birational) 
geometry of $S^{[2]}$.  

For example, a result of Huybrechts and Boucksom 
\cite{Bou,H2} implies that, if the Mori cone of $S^{[2]}$ is closed, then the boundaries are generated by classes of rational curves. Precise numerical and geometric properties of the rational 
curves that are extremal in the Mori cone have been conjectured by Hassett and Tschinkel \cite{HT}.

Uniruled subvarieties of $S^{[2]}$ are important in several aspects: The presence of 
a $\PP^2 \subset  S^{[2]}$ gives rise to a birational map 
(the so-called} Mukai flop,  {\rm cf. \cite{Muk}) to another hyperk\"ahler fourfold and
all birational maps between hyperk\"ahler fourfolds factor through a sequence of Mukai flops 
(see \cite{BHL,HY,W2,WW}). Moreover, uniruled threefolds in $S^{[2]}$ are central in the study
of the birational
K\"ahler cone of $S^{[2]}$ \cite{H2}.

In the particular case of $K3$ surfaces, the study of families of irreducible curves with hyperelliptic normalizations and the loci the curresponding rational curves cover in $S^{[2]}$ is therefore of particular importance. In \cite{fkp2} we study such families. 

Let now $H$ be a globally generated line bundle on $S$ and denote by $|H|^{hyper}$ the subscheme of $|H|$ parametrizing curves with hyperelliptic normalizations. Then, any component of $|H|^{hyper}$ has dimension $\geq 2$ with equality holding if $H$ has no decompositions into moving classes, e.g. if $\Pic S \iso \ZZ[H]$, by \cite[Lemma 5.1]{fkp2}. 

As for concrete examples of such families on a} general {\rm (in the moduli space) primitively polarized $K3$ surface $(S,H)$, here are the ones that are known to us:

(i) $|H|$ contains a two-dimensional family of irreducible curves of geometric genus $p_g=2$, whose general element is nodal, by the nonemptiness of Severi varieties on $K3$ surfaces as a direct consequence of 
Mumford's theorem on the existence of nodal
rational curves on $K3$ surfaces (cf. \cite[pp.~351-352]{MM} or \cite[pp.~365-367]{BPV}) and standard
results on Severi varieties (cf. \cite[Lemma 2.4 and Theorem 2.6]{Tan}; see also
e.g. \cite{CS,F}). In the particular case $p_a(H)=3$, i.e. when $S$ is a smooth quartic in $\PP^3$, the locus in $\Hilb^2(S)$ covered by the associated rational curves is a $\PP^1$-bundle over $S$, 
by \cite[Example 7.7]{fkp2}.

(ii) $|H|$ contains a two-dimensional family of irreducible nodal curves of geometric genus $p_g=3$ with hyperelliptic normalizations, by \cite[Thm. 1.1]{fkp2}. The locus in $\Hilb^2(S)$ covered by the associated rational curves is birational to a $\PP^1$-bundle, by \cite[Cor. 5.3 and Prop. 3.6(ii)]{fkp2}.

(iii) \cite[Prop. 7.7]{fkp2} If $H^2=2(d^2-1)$ for some integer $d \geq 1$, then $\Hilb^2(S)$ contains a uniruled $3$-fold that is birational to a $\PP^1$-bundle. The fibers give rise to a two-dimensional family of curves in $|H|$ with hyperelliptic normalizations of  arithmetic genus $p_a=p_a(H)=d^2$ and geometric genus $p_g =2d-1$. 

(iv) \cite[Prop. 7.2]{fkp2} If $H^2=2(m^2-9m+19)$ for some integer $m \geq 6$, then $\Hilb^2(S)$ contains a $\PP^2$ and the Severi varieties of rational curves in $|\O_{\PP^2}(n)|$, for any $n \geq 1$, give rise to $(3n-1)$-dimensional subschemes $V(n) \subset  |nH|$ parametrizing irreducible curves with hyperelliptic normalizations of arithmetic genus 
$p_a(n)=p_a(nH)=n^2(m^2-9m+19)$ and geometric genus $p_g(n)  =2n-9$.
}
\end{exa}

We have now seen several examples of families as in \eqref{eq:ass1} of 
dimension $2$ on surfaces with $p_g(S) >0$ (in Examples \ref{exa:p2} and \ref{exa:k3}(iv) with $n=1$,
Examples \ref{exa:k3}(i)-(iii) and  the abelian surfaces in Lemma-Example \ref{lemma:abe-biell}) and on Enriques surfaces (the case mentioned in the last lines of Example \ref{exa:enr}).

At the same time we have seen infinite series of examples of  families as in \eqref{eq:ass1}
of arbitrarily high dimensions $\geq 3$ (Examples \ref{exa:p2} and \ref{exa:k3}(iv) with $n \geq 2$, 
and
Examples \ref{exa:baseiper} and \ref{exa:enr}).

In the next section we will see the difference between those ``small'' and ``big'' families.

\section{A Reider-like result} \label{sec:reider}

Consider the additional assumptions
\begin{equation} \label{eq:ass2} 
\dim V  \geq 
\begin{cases} 
3, & \; \mbox{if} \; p_g(S) >0 \; \mbox{or} \; S \; \mbox{is Enriques;} \\ 
5, & \; \mbox{otherwise}.
\end{cases} 
\end{equation}
The following result is an improvement
of \cite[Prop. 4.2]{fkp2}. In fact, the idea of the proof is essentially the same.

\begin{prop} \label{prop:algdec}
Assume \eqref{eq:ass1} and \eqref{eq:ass2}. 
Then there is a decomposition into two effective, algebraically moving classes
\[
[C] = [D_1] + [D_2]
\]
such that, for general $\xi, \eta \in R_V$ (cf. \eqref{eq:defrv}), each with support at two distinct points of $S$,
there are
effective divisors $D'_1  \sim_{alg} D_1$ and $D'_2 \sim_{alg}
D_2$ such that $\xi \subset  D'_1$ and $\eta \subset  D'_2$ and $[D'_1 +
D'_2] \in \overline{V}$, where $\overline {V}$ is the closure of $V$ in the component of the Hilbert scheme of $S$ containing $V$.
\end{prop}

\begin{proof} 
We must have $\dim R_V=2$ or $3$ by \eqref{eq:ass2}. If $p_g(S) >0$, then $\dim V \geq 3$ by \eqref{eq:ass2}, whence $R_V$ is a surface by Proposition \ref{prop:mumford}. If $S$ is Enriques, then there is an unramified double cover  $f:T \khpil S$
such that $T$ is a smooth $K3$ surface, cf. e.g. \cite[VIII, Lemma 15.1(ii)]{BPV}. Therefore, also
$f^{(2)}: \Sym^2(T) \khpil \Sym^2(S)$ is unramified. Hence, the family of rational curves in
$\Sym^2(S)$ given by $\Phi_V$ and $\psi$ as in \eqref{eq:univ} is pulled back to two copies in $\Sym^2(T)$.
Since $\dim V \geq 3$ by 
\eqref{eq:ass2}, we conclude by Proposition \ref{prop:mumford} that these families only cover a surface in $\Sym^2(T)$. Hence 
$R_V \subset \Sym^2(S)$ is a surface as well.

Therefore, in any case, the assumptions \eqref{eq:ass2} guarantee that, for general $\xi, \eta \in R_V$, the locus of points in $V$ parametrizing curves in $\mathcal{P}_V$ passing through $\xi$ and $\eta$ in $R_V$ is at least one-dimensional. 
For general $\xi, \eta \in R_V$,
let $B=B_{\xi,\eta} \subset {V}$ be a smooth curve (not necessarily complete) 
parametrizing such curves and 
\begin{equation} \label{eq:univ10}
\xymatrix{ 
& \widetilde{\C}_B \ar[dl]_{\pi_B} \ar[dr]_{\tilde{\varphi}_B} \ar[r]^{\gamma_B}  & 
\mathcal{P}_B \ar[d]^{\psi_B} \ar[r]^{\Phi_B} & R_V & \hspace{-1cm} \subset  \Sym^2(S) \\
S  & & B & &
}
\end{equation}
the corresponding restriction of \eqref{eq:univ} over $B$. 
 Let $\overline{B}$ be any smooth compactification of $B$. By Mori's bend-and-break technique, as in \cite[Lemma 1.9]{km} or \cite[Cor. II.5.5]{Kol} (see also \cite[Proposition 2.10]{fkp2} for the precise statement we need), there is an extension of the right hand part of \eqref{eq:univ10}
\[ 
\xymatrix{ 
\overline{\mathcal{P}}_B \ar[d]^{\overline{\psi}_B} \ar[r]^{\overline{\Phi}_B} & R_V \\
\overline{B} &
}
\]
such that, for some $b_0 \in \overline{B} \setminus B$ we have
$(\overline{\Phi}_B)_*(\overline{\psi}_B^{-1}b_0) \sup \Gamma_{\xi} + \Gamma_{\eta}$,
where $\Gamma_{\xi}$ and $\Gamma_{\eta}$ are irreducible rational curves on $R_V$ (possibly coinciding) 
such that $\xi \in \Gamma_{\xi}$ and $\eta \in \Gamma_{\eta}$. Let now
\[I:=\{ (x, \xi) \in S \x R_V \; | \; x \in \Supp (\xi) \} \subset  S \x R_V \]
be the incidence variety with projection morphisms $p: I \khpil S$ and $q: I \khpil R_V$.
Then $\dim I=\dim R_V=2$ or $3$  and $q$ is finite of degree two. Consider the commutative diagram
\begin{equation} \label{eq:diagrammonenuovo} 
\xymatrix{ 
\overline{\mathcal{P}}_B \x_{R_V} I =: \hspace{-0,8cm}  &
\overline{\mathcal{P}}'_B \ar[d]_{q'} \ar[r]_{\overline{\Phi}'_B} & I \ar[d]^q \ar[r]_p & S \\
&\overline{\mathcal{P}}_B \ar[d]_{\overline{\psi}_B} \ar[r]^{\overline{\Phi}_B} & R_V & \\
& \overline{B}, & &
}
\end{equation}
where the square is cartesian.
Define $\overline{\pi}_B:=p \circ \overline{\Phi}'_B$.
Note that for $b \in B$ we have $\overline{\pi}_B({q'}^{-1}\overline{\psi}_B^{-1}b)= \pi_B (\tilde{\varphi}_B^{-1}b)$. In particular, $\overline{\pi}_B$ is dominant and generically one-to-one on the fibers over $B$. Therefore $p$ must also be dominant. 

We have 
\[ (\overline{\Phi}'_B)_*({q'}^{-1}\overline{\psi}_B^{-1}b_0) \sup q^{-1}\Gamma_{\xi} + 
q^{-1}\Gamma_{\eta}. \]
Denoting by $b \in B$ a general point and recalling that $\overline{\pi}_B$ is generically one-to-one on the fibers over $B$, we have
\[ C \sim_{alg} (\overline{\pi}_B)_*({q'}^{-1}\overline{\psi}_B^{-1}b) 
\sim_{alg} (\overline{\pi}_B)_*({q'}^{-1}\overline{\psi}_B^{-1}b_0) 
\sup p_*(q^{-1}\Gamma_{\xi}) + p_*(q^{-1}\Gamma_{\eta}) \sup D_{\xi} +D_{\eta},\]
where 
$D_{\xi}:= p(q^{-1}\Gamma_{\xi})$ and $D_{\eta}:= p(q^{-1}\Gamma_{\eta})$.
Now the curves contracted by $p$ are precisely the curves of type $\{x, x+D\}$, for a point $x \in S$
and a curve $D \subset  S$. Since $S$ is not covered by rational curves, and 
$\xi$ and $\eta$ are general, their support on $S$ does not intersect any of the finitely many rational curves $\gamma$ on $S$ with $\gamma.C \leq C^2$. If $q^{-1}\Gamma_{\xi}$ contained a component of the form $\{x, x+D\}$, then, by definition of $q$, we would have
\[ q^{-1}\Gamma_{\xi} = \{x, x+y\}_{y \in D} \cup \{y, x+y\}_{y \in D}, \]
a union of two irreducible rational curves, each being mapped isomorphically by $p$ to
$\Gamma_{\xi}$. Then $p(\{y, x+y\}_{y \in D})=D \subset  S$ would be an irreducible rational curve intersecting $\Supp \xi$, a contradiction. The same argument works for $q^{-1}\Gamma_{\eta}$. 
Therefore, none of the components
of $q^{-1}\Gamma_{\xi}$ nor $q^{-1}\Gamma_{\eta}$ are  contracted by $p$. We therefore have $D_{\xi} \supset \xi $ and $D_{\eta} \supset \eta$,  viewing
$\xi$ and $\eta$ as length-two subschemes of $S$. (Note that $D_{\xi}$ and
$D_{\eta}$ are not necessarily distinct.) Moreover
\[ C \sim_{alg} (\overline{\pi}_B)_*({q'}^{-1}\overline{\psi}_B^{-1}b) = D_{\xi} +D_{\eta} + E_{\xi,\eta}, \]
with $E_{\xi,\eta} \geq 0$, and by construction, $D_{\xi} + D_{\eta}+ E_{\xi,\eta}$ lies in the border of the family 
$\varphi: \C \khpil V$ of curves on $S$, and as such, $[D_{\xi} + D_{\eta}+ E_{\xi,\eta}] \in \overline{V}$, 
where $\overline {V}$ is the closure of $V$ in the component of the Hilbert scheme
of $S$ containing $V$. 
Moreover, as the number of such effective decompositions of $[C]$ is finite (as $S$ is projective), we can find one decomposition
$ [C] = [D_1] + [D_2]$ holding for general $\xi, \eta \in R_V$.   
Since this construction can be repeated for general
$\xi, \eta \in R_V$ and the set $\{x \in S \; |\; x \in \Supp (\xi) \;\; \mbox{for some} \; \xi \in R_V \}$
is dense in $S$, as the curves parametrized by $V$ cover the whole surface $S$, the obtained classes
$D_1$ and $D_2$ must move in an algebraic system of dimension at least one. 
\end{proof}

We now prove an additional, more precise result:

\begin{lemma} \label{lemma:algdec}
Assume \eqref{eq:ass1} and \eqref{eq:ass2}. Then we can find a decomposition as in 
Proposition \ref{prop:algdec}
satisfying the additional properties 
\begin{itemize}
 \item[(a)] $D_1.D_2 \leq p_a-p_g+2$; and
 \item[(b)] there is a reduced and irreducible component of $D_1'$ (resp. $D'_2$) containing $\xi$ (resp. $\eta$).
\end{itemize}
\end{lemma}

\begin{proof} 
Let $\xi$ and $\eta \in R_V$ be general and $[D'_1 + D'_2] \in \overline{V}$ such that 
$\xi \subset  D'_1$ and $\eta \subset  D'_2$ as in Proposition \ref{prop:algdec}. 
Let $D_{\xi} \sub D'_1$, $D_{\eta} \sub D'_2$, $\Gamma_{\xi} \subset R_V$ and $\Gamma_{\eta} \subset R_V$ be as in the proof of Proposition \ref{prop:algdec}, that is, $D_{\xi}:= p(q^{-1}\Gamma_{\xi})$ and $D_{\eta}:= p(q^{-1}\Gamma_{\eta})$. 

If $q^{-1}\Gamma_{\xi}$ were reducible, it would consist of two rational components, each being mapped isomorphically to $\Gamma_{\xi}$ by $q$. Therefore, $\xi$, viewed as a length-two scheme on $S$, would intersect a rational curve $\gamma \subset S$ satisfying $\gamma.C \leq C^2$. As in the proof of Proposition \ref{prop:algdec},
for $\xi$ general this cannot happen. Hence $q^{-1}\Gamma_{\xi}$ is reduced and irreducible and so is 
$D_{\xi}=p(q^{-1}\Gamma_{\xi})$ as well. Of course the same reasoning also works to show that
$D_{\eta}= p(q^{-1}\Gamma_{\eta})$ is irreducible. This proves (b).

Set $D:= D'_1+D'_2$. We now want to show that there is an effective decomposition $D=D_1+D_2$ with $D_{\xi} \sub D_1$, $D_{\eta} \sub D_2$ and $D_1.D_2 \leq p_a-p_g+2$.

We know that a partial desingularization of $D$, say $\widetilde{D}$, which can be obtained by a succession of blowups $f:\widetilde{S} \khpil S$, is a limit of smooth hyperelliptic curves, as
$[D] \in \overline{V}$ by Proposition \ref{prop:algdec}. Let 
$\widetilde{D}_{\xi}$ and  $\widetilde{D}_{\eta}$ denote the strict transforms of  $D_{\xi}$ and $D_{\eta}$, respectively. We now claim that  there is an effective decomposition 
\begin{equation} \label{eq:dec2}
\widetilde{D}=\widetilde{D}_1+\widetilde{D}_2 \; \mbox{with} \; \widetilde{D}_{\xi} \sub \widetilde{D}_1 ,
\;  \widetilde{D}_{\eta}  \sub \widetilde{D}_2  \; \mbox{and} \;  \widetilde{D}_1.\widetilde{D}_2 \leq 2.
\end{equation}

To show \eqref{eq:dec2}, we first write $\widetilde{D}=\widetilde{D}_{\xi} + \widetilde{D}'$. We have a short exact sequence
\[
 0 \hpil   {\omega_{\widetilde{D}_{\xi}}} \hpil    {\omega_{\widetilde{D}}} \hpil    
\omega_{\widetilde{D}'}(\widetilde{D}_{\xi})  \hpil  0.
\]
Since $H^1(\omega_{\widetilde{D}_{\xi}}) \iso H^1(\omega_{\widetilde{D}})$ by Serre duality, the map
$H^0(\omega_{\widetilde{D}}) \khpil  H^0(\omega_{\widetilde{D}'}(\widetilde{D}_{\xi}))$ is surjective, 
whence $|\omega_{\widetilde{D}'}(\widetilde{D}_{\xi})|$ is not birational on $\widetilde{D}_{\eta}$, since
$|\omega_{\widetilde{D}}|$ is $2:1$ on every nonrational component, as
$\widetilde{D}$ is a limit of hyperelliptic curves.

Let now $\widetilde{D}_1 \sub \widetilde{D}$ be maximal with respect to the properties that
$\widetilde{D}_{\xi} \sub \widetilde{D}_1$, $\widetilde{D}_2:=\widetilde{D}-\widetilde{D}_1 \sup \widetilde{D}_{\eta}$ and $|\omega_{\widetilde{D}_2}(\widetilde{D}_1)|$ is not birational on $\widetilde{D}_{\eta}$. 

If $\widetilde{D}_2= \widetilde{D}_{\eta}$, then $\widetilde{D}_1.\widetilde{D}_2 \leq 2$ by
\cite[Prop. 2.3]{CF}, and \eqref{eq:dec2} is proved.

Now assume that $\widetilde{D}_2 \supsetneqq \widetilde{D}_{\eta}$.

If there is a reduced and irreducible component 
$\widetilde{D}_{22} \sub \widetilde{D}_2-\widetilde{D}_{\eta}$ such that
$\widetilde{D}_{22}.\widetilde{D}_{1} >0$, set
$\widetilde{D}_{21}:= \widetilde{D}_{2}-\widetilde{D}_{22} \sup \widetilde{D}_{\eta}$. Then from
\[
\xymatrix{
 0 \hpil   \omega_{\widetilde{D}_{22}}(\widetilde{D}_1)  \hpil    
\omega_{\widetilde{D}_2}(\widetilde{D}_1) \hpil   
\omega_{\widetilde{D}_{21}}(\widetilde{D}_1+\widetilde{D}_{22})  \hpil  0.
}
\]
and the fact that $h^1(\omega_{\widetilde{D}_{22}}(\widetilde{D}_1))=0$, we see, as above, that
$|\omega_{\widetilde{D}_{21}}(\widetilde{D}_1+\widetilde{D}_{22})|$ is not birational
on $\widetilde{D}_{\eta}$, contradicting
the maximality of $\widetilde{D}_1$.

Therefore $\widetilde{D}_1.\widetilde{D}_2= \widetilde{D}_1.\widetilde{D}_{\eta}$ and since 
$|\omega_{\widetilde{D}_2}(\widetilde{D}_1)|$ is  not birational on $\widetilde{D}_{\eta}$, and \linebreak
$H^0(\omega_{\widetilde{D}_{\eta}}(\widetilde{D}_1)) \sub 
H^0(\omega_{\widetilde{D}_2}(\widetilde{D}_1))$, then 
$|\omega_{\widetilde{D}_{\eta}}(\widetilde{D}_1)|$ is not birational on $\widetilde{D}_{\eta}$
either. It follows, using
\cite[Prop. 2.3]{CF} again, that
\[ \widetilde{D}_1.\widetilde{D}_2=\widetilde{D}_1.\widetilde{D}_{\eta} \leq 2, \]
and \eqref{eq:dec2} is proved.

Now let $E_1, \ldots, E_n$ be the (total transforms of the) exceptional divisors of 
$f:\widetilde{S} \khpil S$, so that $K_{\widetilde{S}}= f^*K_S + \sum E_i$,
$\widetilde{D}_1  = f^*D_1- \sum \alpha_i E_i$ and $\widetilde{D}_2  = f^*D_2- \sum \beta_i E_i$, for $\alpha_i, \beta_i \geq 0$, 
where $D_1 \sup \xi$ and $D_2 \sup \eta$. We  compute
\begin{eqnarray*} 
2p_g-2 &  =   & (\widetilde{D}_1+ \widetilde{D}_2).(\widetilde{D}_1+ \widetilde{D}_2 + K_{\widetilde{S}}) \\
     &  =   & (D_1+D_2).(D_1+D_2+K_S) -2\sum \alpha_i.\beta_i +
               \sum \Big (\alpha_i(1-\alpha_i)+\beta_i(1-\beta_i) \Big ), \\
     & \leq & D.(D+K_S) -2\sum \alpha_i.\beta_i =2p_a-2 -2\sum \alpha_i.\beta_i,
\end{eqnarray*} 
whence $\sum \alpha_i.\beta_i \leq p_a-p_g$. Inserting this into
\[ \widetilde{D}_1.\widetilde{D}_2 =(f^*D_1- \sum \alpha_i E_i).(f^*D_2- \sum \beta_i E_i) =D_1.D_2 -
\sum \alpha_i.\beta_i \]
and using \eqref{eq:dec2}, we obtain the desired result $D_1.D_2 \leq p_a-p_g+2$
\end{proof}

A consequence is the following Reider-like result. Note that when $p_a=p_g$, that is, the family consists of {\it smooth} hyperelliptic curves, we retrieve the results of Reider \cite{re}.

We make the following notation: if $C \subset  S$ an irreducible curve with hyperelliptic normalization and $f: \widetilde{S} \to S$ a birational morphism inducing the normalization 
$\nu: \widetilde{C} \to C$, then we define
\[
W_{[2]}(C) := \Big \{ \xi \subset C_{smooth} \; | \; \xi=f_*(Z) \; \mbox{with} \;
Z  \in \mathfrak{g}^1_2(\widetilde{C}) \Big \} \subset  S^{[2]}.
\]

\begin{thm} \label{thm:reider}
Assume \eqref{eq:ass1} and \eqref{eq:ass2}.
Then there is an effective divisor $D$ on $S$ such that  $h^0(D) \geq 2$, 
$h^0(C-D) \geq 2$, $D^2 \leq (C-D)^2$ and 
\begin{equation} \label{eq:inter}
 0 \leq 2D^2 \stackrel{(i)}{\leq}  D.C \leq D^2+p_a-p_g+2 \stackrel{(ii)}{\leq} 2(p_a-p_g+2), 
\end{equation}
with equalities in (i) or (ii) if and only if $C \eqv 2D$.

Furthermore there is a flat family parametrized by a reduced and irreducible complete subscheme $V_D$ of the component of the
Hilbert scheme of $S$ containing $[D]$ with the following property:
for general 
$[C] \in V$ there is a complete rational curve $V_D(C) \sub V_D$ such that for 
general $\xi \in W_{[2]}(C)$, there is a $[D_{\xi}] \in V_D(C)$ such that $\xi \subset   D_{\xi}$. 
\end{thm}

\begin{proof}
We have, by Proposition \ref{prop:algdec} and Lemma \ref{lemma:algdec}, a decomposition into algebraically moving classes $C \sim_{alg} D_1+D_2$ with $D_1.D_2 \leq p_a-p_g+2$. Without loss of generality we can assume that $D_1^2 \leq D_2^2$, or equivalently $D_1.C \leq D_2.C$. We first show that we can assume that $D_1^2 \geq 0$.

Indeed, if $D_1^2 <0$, then the algebraic system in which $D_1$ moves must have a base component $\Gamma >0$. We can write $D_1 \sim_{alg} D_0 + \Gamma$, where $D_0$ moves in an algebraic system of dimension at least one, without base components. In particular $0 \leq D_0^2 = D_1^2-2D_1.\Gamma + \Gamma^2 < \Gamma^2 -2D_1.\Gamma$, so that $\Gamma^2 > 2D_1.\Gamma$. Moreover, we have $\Gamma.D_2 \geq -\Gamma.D_1$ as $C$ is nef. Hence
\[ D_0.(D_2+\Gamma) =D_1.D_2 -\Gamma.D_2+\Gamma.D_1-\Gamma^2 \leq D_1.D_2 +2\Gamma.D_1-\Gamma^2 < D_1.D_2, \]
and we can substitute $D_1$ with $D_0 \sub D_1$, as clearly $\Supp (\xi) \cap \Gamma = \emptyset$ for general $\xi \in R_V$. Therefore, we can assume $D_1^2 \geq 0$.

Combining with  the Hodge index theorem, we get $(2D_1.C) \cdot D_1^2 \leq C^2 \cdot D_1^2 \leq (D_1.C)^2$, so that $2D_1^2 \leq D_1.C$, with equality if and only if $C \eqv 2D_1$. Moreover, $D_1^2 \leq \frac{1}{2}D_1.C = \frac{1}{2}(D_1^2+D_1.D_2) \leq \frac{1}{2}(D_1^2+p_a-p_g+2)$
 yields $D_1^2 \leq p_a-p_g+2$, again  with equality if and only if $C \eqv 2D_1$. 

Finally, from Proposition \ref{prop:algdec} and the fact that we have at most removed base components of the obtained family, it is clear that there is a reduced and irreducible complete scheme $V_{D_1}$ parametrizing curves algebraically equivalent to $D_1$ with the 
property that for general 
$[C] \in V$ and general $\xi \in W_{[2]}(C)$, there is a $[D_{\xi}] \in V_{D_1}$ such that $\xi \subset   D_{\xi}$. For fixed $C$ this yields a complete curve $V_{D_1}(C) \sub V_{D_1}$
such that all $[D_{\xi}] \in V_{D_1}(C)$ for general $\xi \in W_{[2]}(C)$. This gives a natural rational map $C -\khpil V_{D_1}(C)$ inducing a morphism between the normalizations
$\widetilde{C} \khpil \widetilde{V}_{D_1}(C)$ that is composed with the hyperelliptic double cover $\widetilde{C} \khpil \PP^1$. Hence $V_{D_1}(C)$ admits a surjective map from $\PP^1$ and is therefore rational. 

 If $h^0(D_1)=1$, then the variety parametrizing curves algebraically equivalent to $D_1$ is abelian and therefore cannot contain rational curves, cf. e.g. \cite[Prop. 4.9.5]{bl}. Hence $h^0(D_1) \geq 2$.

Substituting $D_1$ with $D_2$ we also obtain $h^0(D_2) \geq 2$, thus finishing the proof.
\end{proof}

In particular, we have a slight improvement of \cite[Cor. 4.7]{fkp2}:

\begin{cor} \label{cor:nodec}
Assume \eqref{eq:ass1} and in addition that there is no decomposition $C \sim_{alg} C_1+C_2$ such that
$h^0(\O_S(C_i)) \geq 2$ for $i=1,2$. 

Then $\dim V \leq 2$ if $p_g(S) >0$ and $\dim V \leq 4$ otherwise.
\end{cor}

The conditions in Corollary \ref{cor:nodec} are for instance satisfied if $\NS(S) \iso \ZZ[C]$.
Theorem \ref{thm:reider} gives  additional restrictions on the existence of such a family as in 
\eqref{eq:ass1} and \eqref{eq:ass2}.
In particular it shows that when the difference $\delta:=p_a-p_g$ is ``small'', then such a family cannot exist unless there are some quite special divisors on the surface (cf. also \cite[Thm. 1]{fkp}, where we in fact show the nonexistence of curves with hyperelliptic normalizations in the primitive linear system $|H|$
with $\delta \leq \frac{p_a-3}{2}$ on a general primitively polarized $K3$ surface $(S,H)$). Unlike the results of Reider, 
Theorem \ref{thm:reider} cannot be used to say that if $[C] \eqv mL$ for some $m>>0$, then 
families as in \eqref{eq:ass2} do not occur, as $p_a$ grows quadratically with $m$, but it shows that the difference $\delta:=p_a-p_g$ must get bigger as $m$ grows. This was already seen in Examples
\ref{exa:p2} and \ref{exa:k3}(iv) above. 

Note that the families in Examples
\ref{exa:p2} and \ref{exa:k3}(iv) for $n \geq 2$ and in Example \ref{exa:enr}  
satisfy the conditions \eqref{eq:ass1} and \eqref{eq:ass2} and that the conditions in Theorem \ref{thm:reider} are satisfied for
$D=H$ in  Examples
\ref{exa:p2} and \ref{exa:k3}(iv) and for $D=2E_1$ in Example \ref{exa:enr}. 

Also note that there is no divisor (in general) satisfying the conditions 
in Theorem \ref{thm:reider} in the two-dimensional families given in Examples
\ref{exa:p2}-\ref{exa:k3} and Lemma-Example \ref{lemma:abe-biell}. We therefore see that the conditions
\eqref{eq:ass2} cannot, in general, be weakened.

\section{Further remarks} \label{sec:appexa}

We make the following observation 
\begin{lemma} \label{lemma:reider}
The divisor class $D$ in Theorem \ref{thm:reider} can be chosen in such a way that
if $D^2 =0$ (resp. $D^2=1$ and $D.C$ is odd), then the general curve parametrized by $V_D$ is reduced and irreducible and $V_D$ is a base point free, complete 
linear pencil with $D.C$ even (resp. ${V_D}$ is a complete linear pencil with one base point $x$ that is a point of every $[C] \in V$). 
\end{lemma}

\begin{proof}
Assume that the general $[D] \in V_D$ is of the form $D=D'+D''$ with $\xi \subset   D'$. Since the number of effective decompositions of $[D]$ in $\NS(S)$ is finite, we can in fact assume that all such $D'$ and $D''$ are algebraically equivalent. For the same reason we can, possibly after moving base components from $D'$ to $D''$, assume that $D'$ is nef, in particular that ${D'}^2 \geq 0$. 
Since ${D'}^2=D^2-2D.{D''}+{D''}^2$, we then get ${D''}^2 \geq 2D.{D''}-D^2$. Moreover we have ${D''}.(C-D) \geq -{D''}.D$ as $C$ is nef. Hence
\begin{eqnarray*} 
D'.(C-D+{D''}) & = & D.(C-D)-D''.(C-D)+D''.D - {D''}^2 \\
& \leq & D.(C-D)+2D''.D - {D''}^2 \leq D.(C-D) +D^2, 
\end{eqnarray*} 
so that $D'.(C-D+{D''}) \leq D.(C-D)$ 
unless when $D^2=1$, ${D''}^2=2D.{D''}-1$ and ${D''}.(C-D)=-{D''}.D$ (whence ${D''}.C=0$). It follows that ${D'}^2=0$, so that $D'.(C-D') \leq \delta+2$, unless $D^2=1$ and $D.C = \delta+3$, in which case $D'.C=\delta+3$.

This proves that the general curve parametrized by $V_D$ can be taken to be reduced and irreducible, possibly upon changing $D$, 
when $D^2=0$. By Theorem \ref{thm:reider}, $|D|$ is a pencil, and as $D^2=0$, we must have that $|D|$ is the whole component of the Hilbert scheme of $S$ containing $[D]$. Therefore ${V_D}=|D| \iso \PP^1$.
Moreover $D.C$ is the degree of the morphism $\widetilde{C} \khpil {V_D} \iso \PP^1$, which is composed with a $g^1_2$, so it must be even. 

We now treat the case $D^2=1$.

We have seen that if the general curve parametrized by $V_D$ is not reduced and 
irreducible, then ${D'}^2=0$ and $D'.C=D.C=\delta+3$, which is odd by assumption. But then, arguing as in the case $D^2=0$, we get that $D'.C$ must be even, a contradiction.
Finally, as the degree of $\widetilde{C} \khpil \widetilde{V}_{D}(C)$ (where, as above, $\widetilde{V}_{D}(C)$ is the normalization of $V_{D}(C)$) must be even,
the family parametrized by $V_D$ must have base points lying on every curve $C$ parametrized by $V$.
Of course there can only be one base point, as $D^2=1$, and the result follows by blowing up $S$ at this base point and using the result for $D^2=0$.
\end{proof}

An additional application of the Hodge index theorem yields the following, which says that in fact when 
$\delta:=p_a-p_g$ is ``small'' the only cases of families of irreducible curves with hyperelliptic 
normalizations of large dimensions arise in fact from $2:1$ rational maps as explained in the introduction and seen in Example \ref{exa:p2}, 
or from elliptic fibrations on the surface (whence $\kod(S) \leq 1$).

\begin{cor} \label{cor:reider2}
Assume \eqref{eq:ass1} and \eqref{eq:ass2}, set $\delta:=p_a-p_g$ and
assume furthermore that $C^2 >(\delta+3)^2$ or that $C^2 =(\delta+3)^2$ and $C \not \eqv (\delta+3)C_0$ for any divisor $C_0$ (resp. the intersection form on $S$ is even, and either
$C^2 >\frac{(\delta+4)^2}{2}$ or $C^2 = \frac{(\delta+4)^2}{2}$ and 
$C \not \eqv \frac{\delta+4}{2}C_0$ for any divisor $C_0$).

Then there is a linear pencil $|D|$, whose general member is a smooth, irreducible elliptic or hyperelliptic curve such that $D^2=0$, $D.C$ is even and $D.C \leq \delta+2$.

Furthermore, in the case when the general member of $|D|$ is hyperelliptic (which is the case if $\kod(S) =2$), there is a
$2:1$ rational map $S \khpil R$ to a smooth rational surface and the family parametrized by $V$ is the pullback of a family of irreducible, rational curves on $R$.
\end{cor}

\begin{proof}
Let $D$ be as in Theorem \ref{thm:reider} and Lemma \ref{lemma:reider}. If $D^2 >0$, then the Hodge index theorem yields
\[ C^2 \leq \frac{(D^2+\delta+2)^2}{D^2} \leq (\delta+3)^2, \]
with equality implying $C \eqv (\delta+3)D$, contradicting our hypotheses. The same reasoning works if the intersection form on $S$ is even, so that $D^2 \geq 2$.

If $D^2=0$, then by Lemma \ref{lemma:reider} $|D|$ is a base point free complete pencil, and for general 
$[C] \in V$ and general $\xi \in W_{[2]}(C)$, there is a $D_{\xi} \in |D|$ such that $\xi \subset   D_{\xi}$.
Therefore, the induced 
morphism $\widetilde{C} \khpil |D| \iso \PP^1$, where $\widetilde{C}$ denotes the normalization of $C$, of degree $D.C$, is composed with the hyperelliptic double cover $\widetilde{C} \khpil \PP^1$. In particular,
$D.C$ is even.

By construction, the general element of $|D|$ has a partial desingularization that admits a $2:1$ map or $1:1$ map onto a $\PP^1$, as it is a component of a limit of smooth rational curves. 
Since $\kod(S) \geq 0$, it cannot be rational, so it is a smooth, irreducible  elliptic or hyperelliptic curve, as $|D|$ is base point free. In the latter case
$S \khpil \PP(\pi_*(\omega_{S/|D|}))$ is the desired $2:1$ rational map, since it maps all members of $V$ generically $2:1$ onto irreducible, rational curves.
\end{proof}

\begin{exa} \label{exa:rankone}
{\rm Assume \eqref{eq:ass1} and \eqref{eq:ass2} with $\NS(S) \iso \ZZ[H]$ such that
$C \eqv mH$ with $m \geq 2$. We must have
$(m-1)H^2 \leq p_a-p_g+2$ by Theorem \ref{thm:reider},
whence 
\[ C^2 =m^2 H^2 \leq \frac{m^2(p_a-p_g+2)}{m-1}.\]
}
\end{exa}

\vspace{0,3cm}

Note that if one can bound $C^2$ and the difference $\delta=p_a-p_g$ is given, one can obtain dimension 
bounds on 
$V$. This follows from the following general result, which is ``folklore''. We include the proof, pointed out to us by C.~Ciliberto,  for lack of a suitable reference. The bounds on $\dim V$ obtained by combining the next lemma
 and Corollary \ref{cor:reider2} are however probably far from being sharp.

\begin{lemma} \label{lemma:dimbound}
Let $V$ be a reduced
and irreducible scheme with $\dim V >0$ parametrizing a flat family of curves on  a smooth projective surface $S$ with $\kod(S) \geq 0$ of 
constant arithmetic and geometric genera $p_a$ and $p_g \geq 2$, respectively, and algebraic equivalence class $C$. 
Then 
\begin{equation} \label{eq:dimbound}
 \dim V \leq \lfloor p_g-\frac{1}{2}K_S.C \rfloor = \lfloor \frac{1}{2}C^2+1-(p_a-p_g) 
\rfloor \leq p_g.  
\end{equation}

In particular, if $\dim V =p_g$, then either $\kod (S) \leq 1$ and $\dim V  =1$, or 
$\kod (S) =0$.
\end{lemma}

\begin{proof}
Note that the equality in \eqref{eq:dimbound} follows from the adjunction formula.

Denote by $C$ a general curve in the family, $\widetilde{C}$ its normalization and $f:\widetilde{C} \khpil S$ the natural morphism. Then we have a short exact sequence 
\begin{equation} \label{eq:df}
\xymatrix{
 0 \ar[r] & \T_{\widetilde{C}} \ar[r]^{df}  & f^*\T_S \ar[r] & \N_f \ar[r] & 0,
}
\end{equation}
defining the {\it normal sheaf} $\N_f$ to $f$. Let $T \subset \N_f$ be the torsion subsheaf and $\overline{\N}_f:= \N_f/T$, which is locally free on ${\widetilde{C}}$. The sections $H^0(T) \subset H^0(\N_f)$ vanish at the generic point of ${\widetilde{C}}$,
cf. \cite[\S\;6]{AC}, whence the tangent space at $[C]$ of the family in which $C$ moves, which corresponds to infinitesimal deformations of $f$ that do not vanish at the generic point of ${\widetilde{C}}$, maps injectively to $H^0(\overline{\N}_f)$, cf. \cite[(3.53)]{ser}. It follows that $\dim V \leq h^0(\overline{\N}_f)$.

From \eqref{eq:df} we have $\deg \overline{\N}_f= \deg \N_f- \length T= -K_S.C +2p_g-2- \length T$.

If now $h^1(\N'_f)=0$, then by Riemann-Roch
\begin{equation} \label{eq:def1}
h^0(\overline{\N}_f) = \deg \overline{\N}_f + 1-p_g = -K_S.C +p_g-1- \length T \leq -K_S.C +p_g-1. 
\end{equation}

If $h^1(\overline{\N}_f) >0$, then by Clifford's theorem we have
\begin{equation} \label{eq:def2}
h^0(\overline{\N}_f) \leq \frac{1}{2}\deg \overline{\N}_f +1 = \frac{1}{2}(-K_S.C +2p_g-2- \length T)+1 \leq 
-\frac{1}{2}K_S.C +p_g.
\end{equation}

Let $\varphi:S \khpil S_0$ denote the morphism to the minimal model of $S$.

If $K_S.C <0$, then $C$ must be an exceptional curve of $\varphi$, so that $K_S.C=C^2=-1$ and $C$ is a smooth rational curve and 
$\dim V=0$, a contradiction.

If $K_S.C \geq 0$, then \eqref{eq:dimbound} follows from \eqref{eq:def1} and \eqref{eq:def2}.
Furthermore, if $\dim V=p_g$, then $K_S.C=0$, so that by adjunction we have 
$2p_g-2 \leq 2p_a-2=C^2$. Moreover, we must have $C =\varphi^*C_0$ for an irreducible curve $C_0 \subset  S_0$ and $K_{S_0}.C_0=0$.

Since $\dim V>0$ by assumption, we must have
$C_0^2 \geq 0$ and the Hodge index theorem implies
$K_{S_0}^2=0$, whence $\kod (S) \leq 1$. Furthermore, if $\dim V \geq 2$, then $C_0^2 \geq 2$, so that
$K_{S_0} \eqv 0$, whence $\kod(S)=0$.
\end{proof}

We conclude the paper by writing out the results in the two simplest cases 
$\delta=p_a-p_g=0$ and $1$.

\subsection{The case $\delta=p_a-p_g=0$} \label{ss:d0}
Assume \eqref{eq:ass1} and \eqref{eq:ass2} with $p_a=p_g$. 
This means that   $V$ parametrizes a flat family of
smooth irreducible curves. We get an effective divisor $D$ as in
Theorem \ref{thm:reider} and Lemma \ref{lemma:reider} and we now consider the various possibilities occuring. 

We first show that the case $D^2=1$ and $C \eqv 2D$ cannot happen.

Indeed, in this case we would have $\dim V=3$ by Lemma \ref{lemma:dimbound}, and by
\begin{equation}  \label{eq:facileD}
0 \hpil \O_S(D-C) \hpil \O_S(D) \hpil \O_C(D) \hpil 0, 
\end{equation}
and the fact that $D-C \eqv -D$ and $D$ is big and nef, as $C$ is nef, we get $h^0(D)=h^0(\O_C(D))=h^0(g^1_2)=2$, so that $\overline{V}_D=|D|$ is a pencil. It has to have one base point, say $x$. But then $|D|$ cuts out a $g^1_1$ on any curve numerically equivalent to $C$ passing through $x$. But, as 
curves algebraically equivalent to $C$ form a family of dimension at least $\dim V=3$, the family of such curves through $x$ has dimension at least $2$, and the surface is uniruled, a contradiction as
$\kod(S) \geq 0$.

Here is a  list of the other possibilities.

{\bf Case $D^2=0$ and $D.C=2$.} Then $|D|$ is a linear pencil cutting out a $g^1_2$ on {\it every} 
smooth curve numerically equivalent to $C$. As in Corollary \ref{cor:reider2}, either this is an elliptic pencil (whence $\kod(S) \leq 1$), or there is a
$2:1$ rational map $S \khpil R$ to a smooth rational surface and the family parametrized by $V$ is the pullback of a family of irreducible, rational curves on $R$.

{\bf Case $D^2=1$ and $D.C=3$.} By  Lemma \ref{lemma:reider}, $|D|$ is a pencil with one base point $x$ lying on every curve parametrized by $V$, and thus cuts out the $g^1_2$ on every member of $V$, and in fact on any smooth curve numerically equivalent to $C$ passing through $x$. 

As $(C-D)^2 \geq D^2$ by Theorem \ref{thm:reider}, we can only have
\[   
C^2=7, \; 8 \; \mbox{or} \; 9, \; \mbox{with $C \eqv 3D$ if $C^2=9$}, 
\]
the latter by the Hodge index theorem.

By Lemma \ref{lemma:dimbound} we have $\dim V \leq 4$ if $C^2=7$ and $\dim V \leq 5$  if 
$C^2=8$ or $9$.

Blowing up at $x$, we reduce to the case above. As $K_S.D \geq 0$, since $\kod(S) \geq 0$ and $D$ moves,
we must have $p_a(D) \geq 2$, so that we can conclude that 
there is a $2:1$ rational map $S \khpil R$  to a smooth rational surface and the family parametrized by $V$ is the pullback of a family of irreducible, rational curves on $R$.

{\bf Case $D^2=2$ and $C \eqv 2D$.} By Lemma \ref{lemma:dimbound} we have $\dim V \leq 5$. Moreover, by \eqref{eq:facileD} we have $h^0(\O_S(D))=h^0(\O_C(D))=h^0(\O_{C'}(D))$ for any $[C] \in V$ and 
any $C' \eqv C$.

Let $|M|$ be the moving part of $|D|$. Note that either $M.C=4$ or $M.C=2$, and the latter implies $M^2=0$ by the Hodge index theorem, so that we can reduce to the case treated above. We can and will therefore assume that $M.C=4$. By the same reasoning, we can assume that $|M|$ is not composed with a pencil.
By the Hodge index theorem, $M^2 \leq 2$ with equality implying $M \sim D$. 

If $h^0(M)=2$, then, by  Theorem \ref{thm:reider}, $|M|$ is a pencil
 with two base points $x$ and $y$ (possibly infinitely near) 
lying on every curve parametrized by $V$, and thus cuts out the $g^1_2$ on every member of $V$, and in fact on any smooth curve numerically equivalent to $C$ passing through $x$ and $y$. It follows that 
$M^2 \geq 2$, whence $M \sim D$. Blowing up at $x$ and $y$, we find as above that there is a $2:1$ rational map $S \khpil R$  to a smooth rational surface and the family parametrized by $V$ is the pullback of a family of irreducible, rational curves on $R$.

If $h^0(M) \geq 3$, then we must in fact have $h^0(M)=3$ and $M^2=2$ by 
Lemma \ref{lemma:dimbound}, whence again $M \sim D$. We have
$h^0(\O_{C'}(D))=3$ for any $C' \in |C|$. Therefore $|\O_{C'}(D)|$ is a $g^2_4$, 
so that in fact every smooth curve in $|C|$ is hyperelliptic. 
It follows from Lemma \ref{lemma:dimbound} that $|D|$ is base point free, so that it defines a $2:1$ morphism $S \khpil \PP^2$, such that
$V$ is the pullback of a subfamily of the $5$-dimensional family of conics in $\PP^2$, as $C \eqv 2D$.

\vspace{0,3cm}

Of course, this is just Reider's result \cite{re} (under stronger hypotheses) if $C^2 \geq 9$. For $C^2 <9$ the results seem to be new with respect to the existing results in
\cite{se,sv,BFL} in the sense that these papers always assume ampleness or very ampleness of 
$\O_S(C)$. (But of course, we have the additional hypotheses on the dimension of the family and on $\kod(S)$.) 

We hope in any case, that the reader may find this treatment of interest because of the completely different approach than the other papers. In particular, we have obtained a ``vector bundle-free'' Reider-like result.

\subsection{The case $\delta=1$} \label{ss:d1}
Assume \eqref{eq:ass1}  with $p_g=p_a-1$. 

We first note that an immediate consequence of Corollary \ref{cor:reider2} is that if $C^2 > 16$ or $C^2=16$ and $C$ is not $4$-divisible in $\NS(S)$, and there is no pencil $|D|$ such that $D.C=2$
(in which case {\it all} curves numerically equivalent to $C$, smooth or not, would carry a 
$g^1_2$, meaning a $2:1$ finite map onto $\PP^1$), then $\dim V \leq 2$ if $p_g(S)>0$ or $S$ is Enriques 
and $\dim V \leq 4$ otherwise.

Under the assumptions \eqref{eq:ass2} we get an effective divisor $D$ as in
Theorem \ref{thm:reider} and Lemma \ref{lemma:reider}. We now consider the various possibilities occuring. 

As above, the case $D^2=1$ and $C \eqv 2D$ cannot happen. Moreover, the case $D^2=0$ and $D.C=3$ does
not happen by Lemma \ref{lemma:reider}.

The cases $(D^2,D.C)=(0,2)$ and $(1,3)$ can be reduced to the case $\delta=0$, in the sense that 
the curves parametrized by $V$ already have $g^1_2$s.

The new cases are:

{\bf Case $D^2=2$ and $C \eqv 2D$.} If $h^0(\O_C(D))=3$, then we can also reduce to the case $\delta=0$. By Lemma \ref{lemma:dimbound} we have $\dim V \leq 4$. 

{\bf Case $D^2=2$ and $D.C =5$.} By the Hodge index theorem, we have $C^2 \leq 12$ and by  
Lemma \ref{lemma:dimbound} we have $\dim V \leq 6$.

{\bf Case $D^2=3$ and $C \eqv 2D$.} By  
Lemma \ref{lemma:dimbound} we have $\dim V \leq 6$.

%%%%%%%%%%%%%%%%%%%%%%%%%%%%%(BIBLIOGRAPHY)%%%%%%%%%%%%%%%%%%%%%%%%%%%%%%%
%
%
%%%%%%%%%%%%%%%%%%%%%%%%%%%%%%%%%%%%%%%%%%%%%%%%%%%%%%%%%%%%%%%%%%%%%%%

\end{document}